%% file: mainv3.tex
\documentclass[titlepage]{amsart}

\usepackage{amssymb}
\usepackage{amsthm}
\usepackage{bbm}
\usepackage{graphicx}

\usepackage[colorlinks=true]{hyperref}
\usepackage[all]{xypic}

\usepackage[final]{epsfig}
\usepackage{psfrag}
\usepackage{epstopdf}

\usepackage{amsaddr}

\usepackage{tikz}
\usetikzlibrary{decorations.pathreplacing}

\usepackage{pgfplots}
\pgfplotsset{compat=1.18}
\usetikzlibrary{arrows.meta}

\newtheorem{thm}{Theorem}[section]
\newtheorem{theorem}[thm]{Theorem}
\newtheorem{proposition}[thm]{Proposition}
\newtheorem{corollary}[thm]{Corollary}
\newtheorem{lemma}[thm]{Lemma}

\theoremstyle{definition}
\newtheorem{definition}{Definition}[section]

\theoremstyle{observation}

\theoremstyle{definition}

\newcommand{\voided}[1]{}

\newcommand{\N}{\mathbb{N}}

{\bf}{\it}
{\bf}{\it}
\newtheorem*{riemann-mapping-theorem}{Riemann Mapping  Theorem}{\bf}{\it}
{\bf}{\it}
{\bf}{\it}
{\bf}{\it}
{\bf}{\it}
{\bf}{\it}
{\bf}{\it}
{\bf}{\it}
{\bf}{\it}
{\bf}{\it}

\newenvironment{pf*}[1]{\proof[#1]}{\endproof}
\usepackage{euscript}

\newcommand{\beq}{\begin{equation}}
\newcommand{\eeq}{\end{equation}}

\newcommand{\aaa}[1]{{{\mathbf{#1}}}}
\newcommand{\crit}{{{\aaa C}}}

\numberwithin{equation}{section}

\newcommand{\RR}{{\mathbb R}}

\newcommand{\ignore}[1]{{}}

\newcommand{\htop}{h_{\mathrm{top}}}

\newcommand{\occurs}{\sqsubset}

 \title[Logistic Computers]{The computational capacity of the Logistic family}

\author{Cristóbal Rojas}
\address{Instituto de Ingeniería Matemática y Computacional,\\ Pontificia Universidad Católica de Chile.\\ Vicuña Mackenna 4860, San Joaquín, Santiago, Chile.   
\\E-mail:  \texttt{\emph{cristobal.rojas{@}mat.uc.cl}}}

\author{Michael Yampolsky}
\address{Department of Mathematics, University of Toronto,\\ 40 St George Street, Toronto, Ontario, Canada. \\ E-mail:  \texttt{\emph{yampol{@}math.toronto.edu}}}

 \thanks{ C.R was partially supported by Grants ANID/FONDECYT  Regular 1230469 and ANID/Basal National Center for Artificial Intelligence CENIA FB210017.  M.Y. was partially supported by NSERC Discovery grant.}
\begin{document}

\begin{abstract}
Viewing computing machines as symbolic dynamical systems is a
powerful idea which has led to many important results. However, the standard
setting of symbolic Turing dynamics is not well suited for studying questions
of typicality, genericity, or robustness. Such natural questions, for example,
as “what is the computational power of a typical computer program?”, “how
common is poly-time computability”, or “how common is universal computation
among all algorithms?” require a different approach. In this paper we
develop a suitable general framework of viewing families of dynamical systems
as computing devices. We study a specific class of such devices given by the
logistic family, and give the answers to the above questions in this setting.
\end{abstract}

\maketitle

\input{preface}
\input{preliminariesv2}

\input{typical_case_proof}

\input{exceptional_case_proof_section}

\input{universal_implies_chaotic}

\bibliographystyle{alpha}
\bibliography{references}

\end{document}

%% file: preface.tex

\section{Introduction}
It is standard in Computer Science to formalize the notion of an algorithm as a {\it Turing Machine} (TM) \cite{Tur}. The definition of a TM predates todays computers, and is awkward from the point of view of applications. However, it is adequate to formalize any code which can be written in, say, {\sl Python} with only a linear increase in the computational cost. TM's operate on symbolic sequences by modifying the symbols within a fixed block of the sequence, and then translating the sequence left or right by a bounded amount. As was observed by C.~Moore \cite{Moore1,moore}, a TM can be naturally seen as a symbolic dynamical system, which he called a generalized shift.

Interpreting the workings of a TM as symbolic dynamics is a poweful idea, which has opened two popular avenues of research. First, it is natural to investigate the properties of TMs as dynamical systems either from analytic or from computational point of view. As an example of the former, in a recent paper \cite{turcomplete}, the authors demonstrate that for  a  class of TMs the topological entropy is positive. As an instance of the latter, in \cite{jeandel} it is shown that the topological entropy of a TM operating on a single symbolic sequence (that is, one-dimensional one) is always computable.

Secondly, one can ask what types of Turing symbolic dynamics can be realized in a given class of smooth/continuous/analytic dynamical systems. Here it is natural to focus on {\it universal} Turing machines -- that is, TMs which can emulate any other TM, and so have ``maximal'' computational power.\footnote{Think of a {\sl Python} interpreter as a concrete example of an algorithm that can emulate any other algorithm.}
For instance, in \cite{koiran}, it was shown that piecewise-monotone smooth maps of the interval or piecewise-linear maps of the square can emulate a universal TM. In another example, in \cite{eulerturing} it was shown that universal TMs can be embedded into Euler flows.

Since the foundational works of Lyapunov and Poincar\'e, the study of smooth dynamical systems focuses on the questions of typicality, genericity, and robustness of dynamical phenomena. Symbolic dynamics is not readily amenable to such questions, and it is difficult to interpret them when TMs are considered. However, they often arise. Consider, for instance, that the authors of \cite{turcomplete} find that their class of TMs with positive entropy contains universal machines. On the other hand, the analysis carried out in \cite{jeandel}  suggests that universal computation could be naturally transcribed into symbolic dynamics with zero entropy. One is naturally led to ask which is more typical either among all TMs or among the universal ones, as well as how typical universal TMs are among all TMs. One could also question how robust  universal computation is under small perturbations. 

There does not seem to be a convenient way of approaching such problems within the setting of symbolic Turing dynamics, or indeed any other conventional interpretation of an algorithm. One struggles to define a ``typical'' computer program, or  make sense of a ``small perturbation'' of an algorithm.  In this paper, we propose a framework which interprets ``conventional'' dynamical systems as computational devices. This is not entirely dissimilar to the literature on finding examples of universal TMs embedded into smooth dynamical systems. However, our focus is different as we propose to view whole classes of smooth dynamical systems as computers. This allows us, in particular, to discuss genericity of computational properties.

A  loose motivation of this approach would be to think of a finite-dimensional family of smooth dynamical systems as an analogue computational device from which we collect a sequence of discrete readings. The device
can be tuned by adjusting finitely many parameters (think of rotating several tuning knobs). As a familiar toy example from one-dimensional dynamics, consider the one-parameter family of quadratic maps of the unit interval
$$f_a(x)=a x(1-x),\; a\in[0,1].$$
The unique critical orbit $\{f_a^k(0.5)\}_{k\geq 0}$ generates the familiar {\it kneading sequence} $\{s_k\}$ from one-dimensional dynamics \cite{kneading}: we assign $s_k=0$ when $f_a^k(0.5)<0.5$ and $s_k=1$ when $f_a^k(0.5)>0.5$. There is, of course, the third possibility when $f_a^k(0.5)=0.5$. However, this can occur only for countably many values of $a$ and can be excluded in a discussion of typical dynamics. When rotating the knob adjusting the value of $a$, there is  zero chance of landing on such a parameter.

Interpreting this sequence of $0$'s and $1$'s as a computation, we can ask, how computationally powerful a typical quadratic map is. As we will see below, this question easily generalizes both to other dynamical systems and to more sophisticated ``encodings''. And it gives non-trivial and insighful answers even in the simplest case considered above.

There is one important convention that we make.  When asking whether the
language generated by a parameter value is decidable, we allow the algorithm
oracle access to that parameter.  This prevents a noncomputable real
parameter from being used as a hidden source of arbitrary information. It also agrees with the intuition of analogue computation, in which the parameters are set externally (with tuning knobs).


%% file: preliminariesv2.tex

\section{Preliminaries}

 \subsection*{Turing machines}
 A (single-tape) Turing Machine \cite{Tur}  is an abstract computing device which, at any given time, can be in one of finitely many internal states. It also has an infinite {\sl working tape}, divided into {\sl cells}. Finitely many of the cells are filled with symbols from a finite alphabet $\Lambda$. There is also a special symbol $B$ which we interpret as a {\sl blank cell} which fills the rest of the tape.
 
 The machine has a \emph{head} which is at any moment positioned over a single cell of the tape. The head can move along the tape and can read, delete and write symbols on the tape, one at a time, according to a finite list of instructions which describe the action to be performed by the machine.
 
 The instructions can only depend on the current internal state of the machine, and the symbol currently being read. It tells the machine whether or not to modify the symbol, which adjacent symbol to read in the next step, and which internal state it should switch to. Some of the internal states are \emph{halting states}. If the machine  enters a halting state then the computation has ended, and the final content of the tape is taken as its output. Otherwise we say the machine does not halt. 
 
Since the set $\Lambda^*$ is countable, we can identify it with $\mathbb{N}$ by some choice of enumeration. In what follows, we will use such an identification to discuss the action of a Turing Machine on natural numbers, where convenient.
 
 We will focus our discussion on yes-or-no decision problems over the inputs of a Turing Machine. 
 Such questions naturally correspond to identifying membership in subsets of $\mathbb{N}$: given a set $A\subseteq \mathbb{N}$, the associated decision problem asks whether a given input $n$ belongs to $A$. A Turing machine $M$ is said to \emph{decide} a set $A\subseteq \mathbb{N}$ if, for every $n\in\mathbb{N}$, $M(n)$ halts and outputs \emph{yes} when $n\in A$ (we say it \emph{accepts} $n$), and \emph{no} otherwise (we say it \emph{rejects} $n$) . Sets with this property are called \emph{decidable}. A set $A\subseteq\mathbb{N}$ is called \emph{recursively enumerable} (r.e.), or \emph{semi-decidable}, if there exists a Turing machine $M$ such that, for every $n\in\mathbb{N}$, $M(n)$ halts and accepts if $n\in A$, and does not halt otherwise. 
All TMs can be enumerated $\{M_1, M_2,\dots\}$ in such a way that one can algorithmically recover a description of $M_e$ from its index $e$, and vice-versa.\footnote{As a practical example of such an ordering, all finite sentences written in {\sl Python} alphabet (including spaces) can be enumerated lexicographically. Most of such ``programs'' would, of course, never run -- which we interpret as halting.} The index of a TM in such an enumeration is usually called its \emph{Goedel number}. We will fix such an enumeration from now on. When run on a given input (including the blank tape), a TM might or might not halt. Let us denote by \[H=\{e \in \mathbb{N} : M_e \text{ halts} \}\] the set of Goedel numbers of Machines that halt on the empty input.
 
 The problem of deciding if a given index $e$ corresponds to a machine that halts on the empty input is known as the \emph{Halting Problem}, and it was shown to be algorithmically unsolvable by Turing \cite{Tur}. 
 
 As we have already stated, we will consider algorithms that can access an external parameter as an oracle. More precisely, we use the oracle convention of \cite{RojasYampolsky2019}.  An oracle for
$a\in\RR$ is a function $\phi$ such that
\[
 |\phi(q)-a|\leq 2^{-q}
\]
and $\phi(q)$ is dyadic.  A query to $\phi$ of precision $q$ costs $q$ time steps -- one may think of it as reading the first $q$ binary digits off an infinite tape on which a binary representation of $a$ is written.  A set $A\subset \mathbb{N}$ is \emph{polynomial-time decidable relative to $a$} if there is an oracle Turing machine that decides whether or not $n\in {A}$ in time polynomial in the input length, with this query cost, for every oracle name $\phi$ of $a$.


\subsection*{Dynamical Systems as computers}

A natural way to view a dynamical system as computing a language or a set is through the symbolic dynamics generated by the orbit of a distinguished point.

Let \(X\) be a space, \(f:X\to X\) a dynamical system in some family $\mathcal{F}$, and \(x_0\in X\) a marked initial point.  A \emph{measurement partition} is a finite partition
\[
\mathcal P=\{P_0,\ldots,P_{r-1}\}
\]
of \(X\), whose atoms are assigned symbols from an alphabet \(A\). The measurement partition records the orbit of \(x_0\) through the symbolic sequence
\[
s(f,x_0,\mathcal P)_n=a
\quad\Longleftrightarrow\quad
f^n(x_0)\in P_a.
\]

Thus, instead of observing the orbit directly, one observes the infinite symbolic sequence
\[
s(f,x_0,\mathcal P)\in A^\N
\]
induced by the partition. The measurement partition determines what information about the orbit is observable. Since our goal is to study the computational power of the dynamics, rather than that of the observation mechanism, we shall always assume that the point $x_0$ and the measurement partition are computationally simple, especially when proving lower bounds (relative to the bounds we would like to prove).
This prevents computational complexity from being encoded in the measurement process itself.

We will need the following notion. 
\begin{definition}[Word occurrence]
Let \(w\in A^*\) and \(x\in A^\N\). We write
\[
w \occurs x
\]
if there exists \(t\in\N\) such that
\[
x_t x_{t+1}\cdots x_{t+|w|-1}=w.
\]
\end{definition}

In what follows, by a \emph{language}, we will mean a subset of the natural numbers. In order to view a dynamical system as computing a language $E\subset \mathbb{N}$, we need to specify a way to encode inputs.

\begin{definition}[Observer]
  An \emph{observer} is a map $\mathcal{W}: \mathbb{N} \to A^*$.  
\end{definition}
 
The idea is for $\mathcal{W}(n)=w_n\in A^*$ to encode the query: does $n \in E \,?$.  The observer must be chosen once and independently from any particular $f$ in the dynamical family under consideration.  
Again, to ensure that the computational content of the set being computed by the system  genuinely arises from the dynamics, rather than being hidden in the observer, we will always assume that observers are computationally simple. For instance, throughout this paper one may take the map \(n\mapsto w_n\) to be computable in polynomial time. The precise complexity bound may depend on the context; its role is to prevent the observer from encoding nontrivial computations. We call such observers \emph{admissible}.

\begin{definition}[Language computed by an observed dynamical system]Let $\mathcal{W}$ be an admissible observer. The language computed by the observed system $(f,x_0,\mathcal P,\mathcal W)$
is
\[
L_{\mathcal W}(f,x_0,\mathcal P)
:=
\{\,n\in\N : w_n \occurs s(f,x_0,\mathcal P)\,\}.
\]
\end{definition}

In other words, the system accepts the input \(n\) whenever the corresponding query word \(w_n\) appears somewhere in the symbolic dynamics generated by the orbit of \(x_0\).

The framework above naturally leads to the following program. Given a class
\(\mathcal F\) of dynamical systems on a space \(X\) with a marked point $x_0\in X$, a measurement partition
\(\mathcal P\), and an admissible observer \(\mathcal W\), one may study
the collection of languages
\[
\bigl\{
L_{\mathcal W}(f,x_0,\mathcal P)
:
f\in\mathcal F\
\bigr\}.
\]
That is, one asks which languages can be computed by the 
dynamics generated by systems in \(\mathcal F\), thus characterizing the computational capacity of the class. The next natural step is to ask how large is the set of systems that compute maximally complex languages. Universality for instance, may occur only at specially tuned parameters, or it may be visible with positive measure, or even occur generically within an open region. The connection between universality and dynamical properties is also relevant, such as, what dynamical properties prevent universality? In what follows, we answer some of these questions for the family of logistic maps.

\subsection*{Logistic Computers} We now specialize the above framework to the logistic family
\[
f_a(x)=ax(1-x), \qquad a\in[0,4].
\]
We shall consider the measurement partition given by the critical point:
\[
\mathcal P:=\bigl\{[0,0.5),[0.5,1]\bigr\}
\]
so that the symbolic dynamics induced by \((f_a,0.5,\mathcal P)\) is precisely the classical kneading sequence:

\begin{definition}
For \(a\in[0,4]\), the \emph{kneading sequence} \(K(a)\in\{0,1\}^{\N}\) is given by
\[
K(a)_k=
\begin{cases}
0,&f_a^{k+1}(0.5)<0.5,\\[2mm]
1,&f_a^{k+1}(0.5)\ge 0.5.
\end{cases}
\]
\end{definition}

We call \(a\) an \emph{excluded parameter} if
\[
f_a^k(0.5)=0.5
\]
for some \(k\in\N\). There are only countably many such parameters.
Given an admissible observer $\mathcal{W}$,
the language computed by the logistic map \(f_a\) is therefore
\[
S_{\mathcal W}(a)
:=
L_{\mathcal W}(f_a,0.5,\mathcal P)
=
\{n\in\N:\; \mathcal{W}(n)\occurs K(a)\}.
\]

Equivalently, the input \(n\) is accepted precisely when the query word \(\mathcal{W}(n)=w_n\) appears in the kneading sequence of \(f_a\). We shall refer to the pair \((a,\mathcal W)\) as a \emph{logistic computing device}.

\section{Main results}

Our main result is a sharp typical-vs-exceptional classification of logistic computing devices. On the one hand,
for any admissible observer \(\mathcal W\), the language \(S_{\mathcal W}(a)\)
is polynomial-time decidable for Lebesgue-almost every parameter.
On the other hand, the exceptional parameters are
computationally maximal: there exists a poly-time computable observer \((\mathcal W)\) such that every r.e.\ set is realized \emph{exactly} as \(S_{\mathcal W}(a)\)
for some computable parameter \(a\), and the parameter is produced effectively from the G{\"o}del index of the
r.e.\ set. This gives a precise sense in which Turing-completness is possible but atypical in low-dimensional
dynamics: it occurs on exceptional parameters, while almost every parameter is
algorithmically tame.

We further show that Turing-complete parameters must be \emph{chaotic}, in the sense that they necessary lie in the positive-topological-entropy regime, formally linking Turing universality to a standard notion of dynamical complexity. 
We now state the main results precisely. Firstly, we give a formal  definition of Turing completeness.

\begin{definition}[Turing complete parameters]
  Let us say that a set $A\subset\N$ is {\it many-one reducible} from a set $B\subset N$ if there exists a computable function $f$ on natural numbers such that $n\in A$ if and only if $f(n)\in B$. Intuitively, this means that $A$ is ``easier'' to decide than $B$, as there is a computational recipe for determining membership in $A$ if membership in $B$ is known.
 A r.e. set is {\it r.e.-complete} under many-one reductions if any other r.e. set is reducible from it. We will say that a logistic computing device $(a,\mathcal W)$ is {\it Turing-complete}  
if \(S_{\mathcal W}(a)\) is r.e.-complete.
\end{definition}

\begin{theorem}[\bf Main classification theorem]\label{thm:main-classification}
Under the definitions given above:
\begin{itemize}
  \item[(i)]  For Lebesgue-almost every parameter \(a\in[0,4]\) and every poly-time computable observer \(\mathcal W\), the language \(S_{\mathcal W}(a)\subseteq \N\) is
  polynomial time decidable relative to $a$.  

  \item[(ii)]  There exists a poly-time computable observer \(\mathcal W\) and a computable parameter $a_U$ such that the logistic computing device $(a_U,\mathcal W)$ is {\it Turing complete}. 
  
\end{itemize}
\end{theorem}

We note that the positive part of our result is
nonuniform in $a$: the oracle machine and its polynomial bound may depend on
the fixed parameter. 

The next observation shows that, within the logistic family, only chaotic systems can be Turing complete.

\begin{theorem}[\bf Universality implies  chaos]\label{thm:entropy-universality}
Let $\mathcal W$ be an admissible observer. If the logistic computing device \((a,\mathcal W)\) is Turing complete, then 
\(
\htop(f_a)>0.
\)
\end{theorem}

\section{Related work}

Since Moore's interpretation of Turing machines as generalized shifts
\cite{Moore1,moore}, there have been many works dealing with the question of which computational processes can be embedded into more classical dynamical systems. For example, Koiran showed that universal computation can be realized by
piecewise-monotone smooth interval maps and by piecewise-linear maps of the square \cite{koiran}.  Gra{\c c}a and Zhong proved that one-dimensional analytic
maps, and two-dimensional analytic ordinary differential equations, suffice
for robust Turing simulation \cite{GracaZhong2023} over unbounded domains. In fluid dynamics,
Cardona et al. proved Turing
completeness results for Euler flows \cite{CardonaMirandaPeraltaSalas2022,
CardonaMirandaPeraltaSalasPresas2023}.  More recently, Cotler and Rezchikov
formulated a general complexity-theoretic framework for computational
dynamical systems, emphasizing that a simulation must specify not only the
dynamics but also low-complexity encoders and decoders
\cite{CotlerRezchikov2024}.

The observer formalism used in this paper is a symbolic instance of this
general viewpoint.  A class of systems, a measurement partition, and an
observer specify how computation is encoded and read from an orbit.  The
distinctive feature of our results is that these choices are rigid: the
measurement is the natural critical itinerary, the coding \(n\mapsto w_n\) is
fixed once and for all and polynomial-time computable, and the readout is
only the occurrence of \(w_n\) in the kneading sequence.  Thus, besides asking whether a system is capable of simulating a Turing machine, we asks
how often such computational power appears inside a natural parameter
family.

On the dynamics side, kneading theory gives the basic symbolic
description of unimodal dynamics: the itinerary of the critical point
determines the admissible symbolic dynamics of the map
\cite{kneading}; see also \cite{deMeloVanStrien1993}.  This
naturally leads to the study of unimodal languages as formal languages.
Grassberger and Hao studied forbidden words and grammatical structure for
one-humped maps \cite{Grassberger1988,Hao1989}.  K{\r u}rka gave a systematic
classification of the possible language complexities of unimodal systems,
including nonrecursive languages, recursive non-context-sensitive languages,
and context-sensitive nonregular languages \cite{Kurka1996}.  Further
refinements for eventually periodic and aperiodic kneading sequences were
obtained by Wang, Xie, and Yang \cite{WangXie1994,WangYangXie1999}.

The computability issue depends strongly on what information is available
about the parameter. Earlier  noncomputability results were  usually formulated
in an absolute sense: the parameter is fixed, but not provided as an oracle.
Then a noncomputable parameter can itself be the source of noncomputability
in the symbolic dynamics.  This is the viewpoint in Friedman's work on
one-dimensional families, where  noncomputable behavior occurs on a set of positive measure
\cite{Friedman1991}.  Similarly, in K{\r u}rka's nonrecursive examples one
chooses admissible kneading data of the desired complexity; without explicitly separating this complexity from that of the selected parameter. 

Our formulation factors out this external source of information by studying
decidability relative to the parameter.  What remains is a dynamical question: can the critical orbit
of \(f_a\), once the parameter is given, still generate universal
 behavior through the fixed (and simple) observer?

The answer separates typical and exceptional behavior.  For almost every parameter, the observer language is efficiently decidable relative to the
parameter. Genuinely dynamical universality is, on the other hand, a zero-measure phenomenon.

%% file: typical_case_proof.tex
\section{The typical case: proof of Theorem~\ref{thm:main-classification}(i)}

Let us recall that a logistic map \(f_a\) is called \emph{regular} if it has an attracting periodic cycle in \([0,1]\). In this case, the attracting cycle is unique. The attracting basin contains the critical point \(c=1/2\) as well as all of the other points in $[0,1]$, except finitely many. We shall use the following theorem of Lyubich \cite{DynQuad,RegStoch}:

\begin{theorem}\label{th-attractor}
For Lebesgue-almost every \(a\in[0,4]\), one of the following two possibilities holds:
\begin{enumerate}
    \item \(f_a\) is regular;
    \item there exists an interval \(I_0\ni c\) and an integer \(k\geq 1\) such that
    \[
    f_a^k(I_0)=I_0,
    \]
    the intervals
    \[
    I_j:=f_a^j(I_0), \qquad 0\leq j<k,
    \]
    are pairwise disjoint, and the critical orbit
    \[
    \{f_a^m(c):m\geq 0\}
    \]
    is dense in
    \[
    \bigcup_{j=0}^{k-1} I_j.
    \]
\end{enumerate}
Moreover, each alternative occurs for a set of parameters of positive measure.
\end{theorem}
 
In order to show polynomial time decidability for a.e. parameter, we will use the following result by Avila et al. which establishes polynomial recurrence of the critical point (see \cite{AvilaMoreira2005}). This is what will allow us later to read the kneading sequence in polynomial time. 

\begin{theorem}\label{thm:polynomial_recurrence}
    Almost every non-regular parameter has polynomial recurrence of the critical orbit with exponent 1: 

    \begin{equation*}\label{eq:poly-recurrence}
 \limsup_{t\to\infty}
 \frac{-\log|f_a^t(c)-c|}{\log t}=1
\end{equation*}

In other words, the set of $t$ such that $|f_a^t(c)-c|<t^{-\beta}$ is finite whenever $\beta>1$. 
\end{theorem}

We will use the previous result in the form of the following direct corollary:

\begin{corollary}\label{cor:polytime_kneading}
For almost every non-regular parameter $a$, there exists an oracle algorithm  which, on input $N$, computes
\[
 K(a)_0K(a)_1\cdots K(a)_{N-1}
\]
in time polynomial in $N$.
\end{corollary}

\begin{proof}
    
For a non-regular parameter $a$ satisfying Theorem \ref{thm:polynomial_recurrence}, there exists a constant $C_a$ such that, for every time $t$ and $t\leq N$:
\begin{equation*}\label{eq:recurrence-lower}
 |f_a^t(c)-c|\geq C_a N^{-2}.
\end{equation*}
Assume now we perform every arithmetic operation with $q$ bits of precision, and denote by $e_t$ the maximal error made when computing $f^t_a(c)$. Using a uniform modulus of continuity of $f_a$ one shows that:

\begin{equation*}\label{eq:error-propagation}
 e_t\leq C_0 4^t2^{-q}
\end{equation*}
for some absolute constant $C_0$.

Therefore, in order to determine whether $f_a^t(c)>c$ or $f_a^t(c)<c$ (and thereby the $t$-th symbol in $K(a)$), it is enough to choose $q$ so that

\[
 C_0 4^N2^{-q}<\frac12 C_aN^{-2}.
\]
For example,
\begin{equation*}\label{eq:precision-bound}
 q\geq 2N+2\log_2N+\log_2(2C_0/C_a)
\end{equation*}
is sufficient.  Thus only $O_a(N)$ bits of the parameter and of working
precision are required, and the involved operations can clearly be performed in polynomial time.  
\end{proof}

We also note (cf. \cite{YRAdvances19}):
\begin{proposition}\label{comp-attractor}
In case (i) of Theorem \ref{th-attractor} the attracting periodic orbit is uniformly computable given $a$ with a non-uniform polynomial time complexity bound. In case (ii) the boundary points of the interval $I_0$ are computable given $a$. The latter computation is uniform given $a$, $k$, and is similarly, non-uniformly poly-time.
\end{proposition}

We will also need the following simple observation.

\begin{lemma}\label{lem:decidable-subword-implies-decidable-language}
Let \(x\in\{0,1\}^{\mathbb N}\) be an infinite binary sequence whose finite subword language
\[
\mathcal L(x):=\{w\in\{0,1\}^*: w\occurs x\}
\]
is polynomial time decidable. Then, for every polynomial time computable observer
\[
\mathcal W=(w_0,w_1,w_2,\ldots),
\]
the set
\[
S_{\mathcal W}(x):=\{n\in\mathbb N : w_n\occurs x\}
\]
is polynomial time decidable.
\end{lemma}

\begin{proof}
Given \(n\), compute \(w_n\). Since \(\mathcal L(x)\) is poly-time decidable, we can decide whether
\[
w_n\in\mathcal L(x)
\]
in polynomial time, which is equivalent to \(n\in S_{\mathcal W}(x)\).
\end{proof}


\subsection*{The regular case}

Suppose first that \(f_a\) is regular. Since the attracting periodic orbit attracts the critical point \(c\), the itinerary of \(c\) with respect to the critical partition
\[
[0,c)\cup[c,1]
\]
is eventually periodic. Hence the kneading sequence \(K(a)\) is eventually periodic and therefore $\mathcal{L}(K(a)$ is poly-time decidable. 

\ignore{
\begin{lemma}\label{lem:eventually-periodic-subword-decidable}
If \(x\in\{0,1\}^{\mathbb N}\) is eventually periodic, then its finite subword language \(\mathcal L(x)\) is decidable.
\end{lemma}

\begin{proof}
Write \(x=uvvvv\cdots\) for some finite words \(u\in\{0,1\}^*\) and \(v\in\{0,1\}^+\). Given a word \(w\), let \(\ell=|w|\). If \(w\) occurs in \(x\), then it occurs starting either inside the finite prefix \(u\), or inside the periodic tail. It is therefore enough to check whether $w$ occurs in the finite word \(u\,v^{N}\) 
for \(N\) large enough.  This is decidable by finite search.
\end{proof}
Combining Lemmas~\ref{lem:decidable-subword-implies-decidable-language} and~\ref{lem:eventually-periodic-subword-decidable}, we obtain that for every computable observer \(\mathcal W\), the language \(S_{\mathcal W}(a)\) is decidable whenever \(f_a\) is regular.
}

\subsection{The non regular case}

We prove the following:

\begin{proposition}
Let \(a\) be a non regular parameter satisfying the second possibility in Theorem \ref{th-attractor} and the conclusion of Corollary \ref{cor:polytime_kneading}. Then \(\mathcal L(K(a))\) is poly-time decidable relative to
\(a\).  Hence, for every poly-time computable observer \(\mathcal W=(w_n)\), the set
\[
        S_{\mathcal W}(a)
        =
        \{n:\ w_n\in \mathcal L(K(a))\}
\]
is poly-time decidable relative to \(a\).
\end{proposition}

\begin{proof}
The attractor is a cycle of intervals
\[
        A_a=I_0\cup I_1\cup\cdots\cup I_{p-1},
        \qquad f_a(I_j)=I_{j+1 \!\!\!\pmod p},
\]
with \(0.5\in I_0\), and the critical orbit is dense in \(A_a\).  We take the minimal such
\(p\).  By  Proposition \ref{comp-attractor},  the endpoints of the intervals in the cycle are poly-time computable relative to
\(a\).

For \(1\leq j<p\), the interval \(I_j\) lies on one side of \(0.5\).  Let
\(b_j\) be this side symbol and set \(B=b_1\cdots b_{p-1}\).  The word \(B\)
is computable relative to \(a\).  If \(g=f_a^p|_{I_0}\) and \(r_m\) is the
side of \(g^m(0.5)\), then the kneading sequence has the form
\[
        K(a)=B r_1 B r_2 B r_3\cdots,
\]
and is therefore enough to decide the language of the return sequence
\(R=r_1r_2\cdots\). We will need the following Lemma:

\begin{lemma}\label{lemma:admissibility_condition}
    For binary words of the same length, let $\preceq$ be the parity-lexicographic order: at the first disagreement, if the common prefix has an even number of 1's, use $0<1$ and reverse the order otherwise. Let $I_0=[l,r]=[g^2(0.5),g(0.5)]$. Then a word $u=u_0\dots u_{m-1}$ is realized as the itinerary under $g$ of a nonempty open interval of points in $I_0$ if and only if 
    \begin{equation}\label{condition}
          \sigma K_g(a)|_{m-i}\preceq u_i\dots u_{m-1}\preceq K_g(a)|_{m-i} \qquad \text{ for all } \quad 0\leq i < m. 
     \end{equation}
    where $\sigma$ denotes the operation that deletes the first symbol and $K_g(a)$ denotes the itinerary under $g$ of $0.5$.  
\end{lemma}
\begin{proof}
We assume $g$ has a maximum at $0.5$. The other case can be handled similarly. 
Let the itinerary map $\iota_g(x)$ of a point $x \in I_0$ be defined by
\[
\iota_g(x)_j=\mathbf1_{\{g^j(x)\ge 0.5\}},
\qquad K_g=\iota_g(g(0.5)).
\] 

It is nondecreasing with respect to the
parity-lexicographic order \cite[Chapter~II, Proposition~3.1(a)]{deMeloVanStrien1993}. Consequently, the same holds for finite itinerary prefixes. If $u$ is realized by the itinerary of $x\in[l,r]$, then (\ref{condition}) must hold since $\iota_g(l)=\sigma K_g(a)$ and $\iota_g (r)=K_g(a)$. Suppose now that $u$ satisfies (\ref{condition}). We use induction on $m$. The case $m=1$ is trivial since both $0$ and $1$ are always realized by some $l<x<r$. Assume it holds for $m$ and let $u=bv$ with $|v|=m$, so $v$ is realized by points in a nonempty open interval $C(v)$ inside $I_0$. If $b=1$, since the decreasing branch of $g$ maps $(0.5,r)$ onto $(l,r)$, the preimage of $C(v)$ along that branch realizes $u$ on an nonempty open interval. Suppose now $b=0$ and set $z=g(l)$ so that the increasing brach maps $(l,0.5)$ onto $(z,r)$. Since $\iota_g(l)=0\iota_g(z)$, by (\ref{condition}) we have that
\[
\iota_g(z)|_{m}\preceq v. 
\]
If the inequality is strict, then $C(v)$ lies entirely to the right of $z$. If equality holds, since the first $m$ iterates of $z$ avoid $0.5$, the same itinerary is realized by points in a sufficiently small neighborhood to the right of $z$. In either case, $C(v)\cap(z,r)$ contains a nonempty open interval, and its inverse image under the increasing brach realizes $u$. 
\end{proof}

To finish the proof, we note that since $0.5$ is dense in $I_0$ under $g$, a word $u$ is realized as the itinerary of a point if and only if $u$ appears in $K_g(a)$. By Lemma \ref{lemma:admissibility_condition}, this can be decided from $K_g(a)|_{m+1}$ with $O(m^2)$ symbol comparisons. Since $K_g(a)$ is polytime computable relative to $a$, it follows that both $\mathcal{L}(K_g(a)$ and $\mathcal{L}(K(a))$ are decidable in polynomial time relative to $a$. Applying this decision procedure to the 
computable word \(w_n\) proves the assertion for every computable observer.
\end{proof}

%% file: exceptional_case_proof_section.tex
\section{Contructing universal parameters: proof of Theorem \ref{thm:main-classification} (ii)}

\begin{figure}[ht]
\centering
\input{logistic_orbit_adds_ten_tikz}
\caption{
  The figure shows \(f_a\), $a=a_*$, the lightly drawn graph of
\(g_a=f_a^3\), the critical \(g_a\)-orbit
\(0.5=c\mapsto\beta_a'\mapsto\beta_a\), and two three-step orbit portraits
starting from representative points in \(L_a\) and \(R_a\).  These orbits show
that both branches first move to the right of \(c\), then to the left of
\(c\), before returning to the horseshoe under \(g_a\).  Thus passing from
the \(g_a\)-itinerary to the ordinary kneading sequence replaces each
\(g_a\)-symbol by the corresponding block \(\Theta(0)=010\) or
\(\Theta(1)=110\).}
\label{fig:theta-adds-ten}
\end{figure}
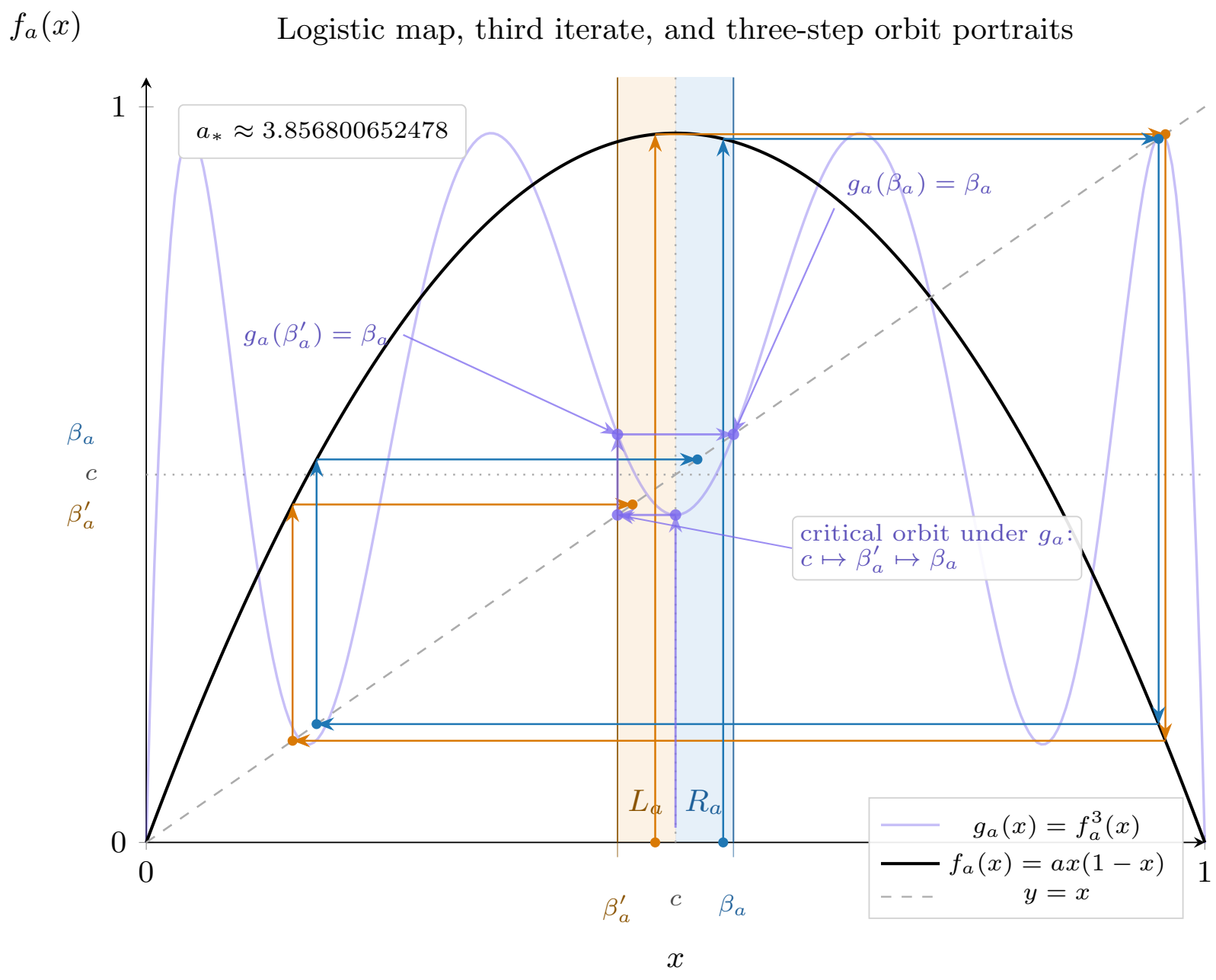

For $a\in(0,4]$, consider the third iterate 
\[
        g_a\equiv f_a^3.
\]
We now use the following convention:
\[
        \iota(x)=
        \begin{cases}
        0,& x<c,\\
        1,& x\geq c
        \end{cases}
\]
so that $K(a)_j=\iota(f_a^{j+1}(c))$. We start by noting that there exists a parameter value \(a_\ast \in (3.85,4)\) such that the following holds: 
\[
        g_{a_\ast}(0.5)\neq g_{a_\ast}^2(0.5)=g_{a_\ast}^3(0.5).
\]
If we set 
\[
        \beta_\ast=g_{a_\ast}^2(0.5),
        \qquad
        \beta_\ast'=g_{a_\ast}(0.5),
\]
then
\[
        \beta_\ast'<c<\beta_\ast,\qquad
        \beta_\ast'=1-\beta_\ast,\qquad
        g_{a_\ast}(\beta_\ast)=\beta_\ast,
\]
and denoting $I_{a_\ast}\equiv [\beta'_\ast,\beta_\ast]\ni 0.5$, we have
$$g_c(I_\ast)=I_\ast,$$
where both endpoints of $I_\ast$ map to $\beta_\ast$. The restriction $g_\ast|_{I_\ast}$ maps both halves $L_\ast=[\beta_\ast',0.5]$ and $R_\ast=[0.5,\beta_\ast]$  of the  interval $I_\ast$ onto the whole $I_\ast$ in a monotone fashion (that is, it folds $I_\ast$ over itself). See Figure \ref{fig:theta-adds-ten} for an illustration.

The fixed point \(\beta_\ast\) of \(g_{a_\ast}\) is repelling.  Hence it has a
real-analytic continuation \(\beta_a\) for \(a\) close to \(a_\ast\) and we again set 
\[
        \beta_a'=1-\beta_a,
        \qquad
        I_a=[\beta_a',\beta_a].
\]

Choose a one-sided interval \((a_\ast,a_\ast+\varepsilon_0)\), with
\(\varepsilon_0>0\) small, so that for every \(a_\ast<a<a_\ast+\epsilon_0\) the following holds:
\begin{enumerate}
    \item The equation
\[
        g_a(x)=\beta_a'
\]
has exactly two solutions \(l_a,r_a\) in \(I_a\), with $\beta_a'<l_a<c<r_a<\beta_a$.
\item \(g_a\) is monotone on
\[
        L_a=[\beta_a',l_a],
        \qquad
        R_a=[r_a,\beta_a],
\]
and it maps $L_a$ and $R_a$ both onto $I_a$.
\item For every $x\in L_a\cup R_a$ we have
\[
        f_a(x)>c,\qquad f_a^2(x)<c.
\]
\end{enumerate}

Therefore a \(g_a\)-symbol determines three ordinary kneading symbols:
\[
        x\in L_a
        \quad\Longrightarrow\quad
        \iota(x)\iota(f_a(x))\iota(f_a^2(x))=010,
\]
whereas
\[
        x\in R_a
        \quad\Longrightarrow\quad
        \iota(x)\iota(f_a(x))\iota(f_a^2(x))=110.
\]

We therefore introduce the block code
\[
        \Theta(0)=010,\qquad \Theta(1)=110.
\]

\begin{figure}[ht]
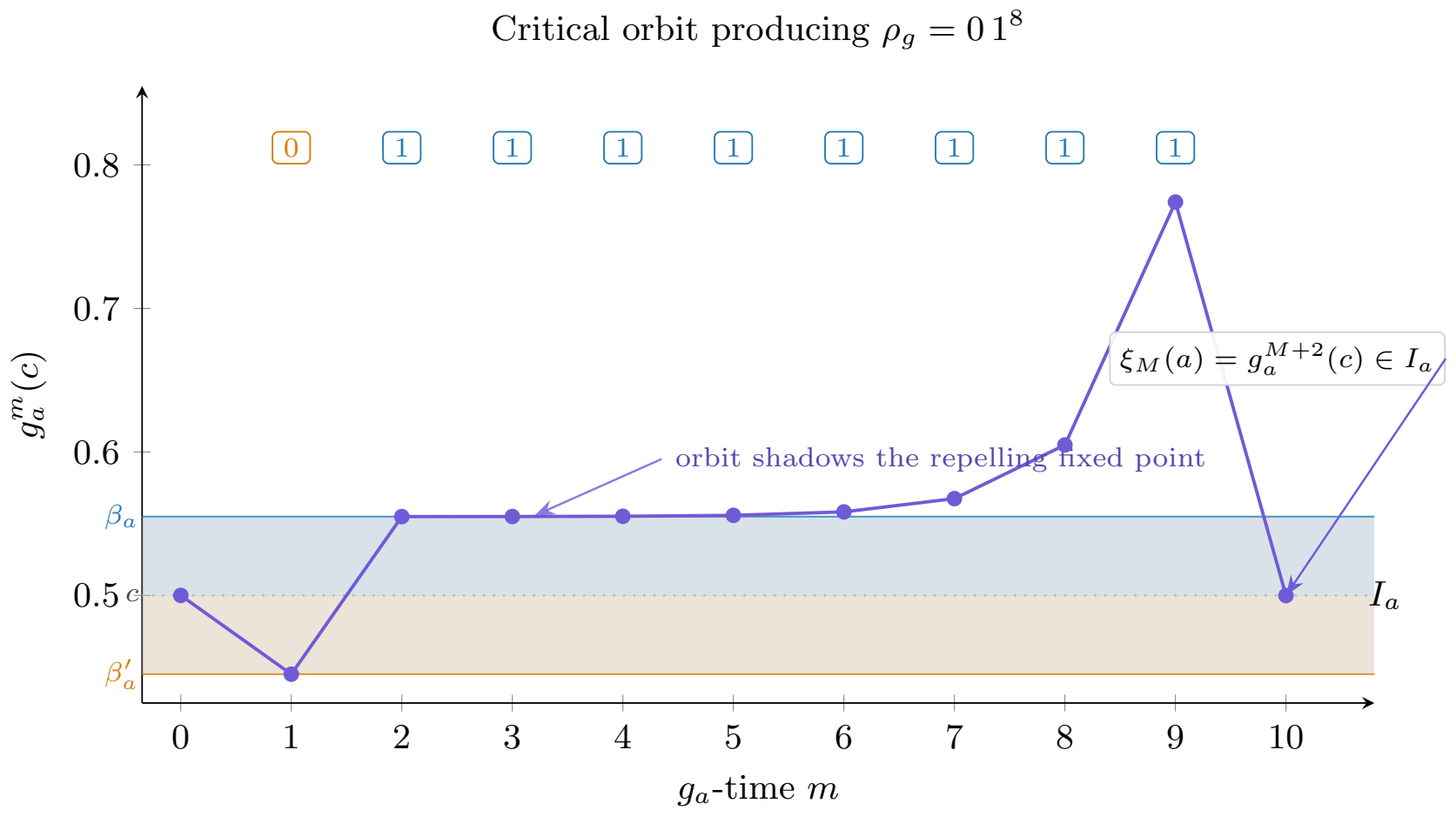

\centering
\include{logistic_prefix_landing_tikz}
\caption{The critical \(g_a=f_a^3\)-orbit in an illustrative landing window
with \(M=8\).  First \(g_a(c)\) lies on the left of \(c\), giving the symbol
\(0\).  Then the orbit lands very near the repelling fixed point
\(\beta_a\), remains to the right of \(c\) for \(M\) iterates, and finally
lands back in \(I_a\) at \(\xi_M(a)=g_a^{M+2}(c)\).  This is the mechanism
which gives the prefix \(\rho_g=0\,1^M\), and
\(\rho_f=10\,010\,(110)^{M-1}111\) in the ordinary kneading sequence.}
\label{fig:critical-prefix-landing}
\end{figure}

The maximal invariant set
\[
        C_a=\{x\in I_a:\ g_a^n(x)\in L_a\cup R_a
        \text{ for every }n\geq0\}
\]
is a Cantor set.  Since each of the two branches \(L_a,R_a\) maps
monotonically onto \(I_a\), the itinerary map
\[
        \pi_a:C_a\longrightarrow\{0,1\}^{\mathbb N},
        \qquad
        \pi_a(x)_n=
        \begin{cases}
        0,& g_a^n(x)\in L_a,\\
        1,& g_a^n(x)\in R_a,
        \end{cases}
\]
conjugates \(g_a|_{C_a}\) to the full one-sided shift.  In particular, for every finite word
\(u\in\{0,1\}^m\), the cylinder
\[
        I_a(u)=\{x\in I_a:\pi_a(x)_0\cdots \pi_a(x)_{m-1}=u\}
\]
is a compact interval, and \(g_a^m\) maps \(I_a(u)\) monotonically onto
\(I_a\). Now, we note that the \(g_{a_\ast}\)-itinerary of the critical orbit is
\[
        0\,1\,1\,1\dots,
\]
because $g_{a_\ast}(0.5)=\beta_\ast'<0.5$, $g_{a_\ast}^2(0.5)=\beta_\ast>0.5$ and $g_{a_\ast}(\beta_\ast)=\beta_\ast.$ We now choose a finite prefix of the form $\rho_g=0\,1^M$ for some large $M$ and let 
\[
 \xi_M(a)=g_a^{M+2}(0.5).
\]
By the transversality of the critical orbit, for all
sufficiently large \(M\), the map \(\xi_M(a)\) crosses the moving interval
\(I_a=[\beta_a',\beta_a]\) monotonically on a small parameter
window \(U\subset (a_\ast, a_\ast+\epsilon_0)\). For these parameters, the
orbit segment
\[
        g_a^2(0.5),\,g_a^3(0.5),\,\ldots,\,g_a^{M+1}(0.5)
\]
stays to the right of \(0.5\), while the next point
\(\xi_M(a)=g_a^{M+2}(0.5)\) sweeps across \(I_a\) (See Figure \ref{fig:critical-prefix-landing} for an illustration).  This is why the
prefix $\rho_g$ can be fixed as \(0\,1^M\). In the
kneading sequence of the original map \(f_a\), this fixed prefix can be chosen to be
\[
        \rho_f
              =10\,010\,(110)^{M-1}111.
\]
The initial \(10\) is forced by \(f_a(0.5)>0.5\) and \(f_a^2(0.5)<0.5\). The three 1's at the end will play a role in the proof below. They can be enforced by choosing the descending return branch on which, for
\[
x=g_a^{M+1}(0.5),
\]one has $
x>c, \,\,f_a(x)>c, \,\, f_a^2(x)>c,$ throughout \(J\). Concretely, there are parameters
\[
        a_M^-<a_M^+,\qquad a_M^-,a_M^+\in U,
\]
such that \(\xi_M(a_M^-)=\beta_a\), \(\xi_M(a_M^+)=\beta_a'\), and  for every \(a\in J=[a_M^-,a_M^+]\)
we have
\[
        g_a(0.5)<0.5,\qquad
        g_a^j(0.5)>0.5\quad \text{for }2\leq j\leq M+1,
        \qquad
        \xi_M(a)\in I_a.
\]
Thus every \(a\in J\) has the same finite \(g_a\)-itinerary prefix
\[
        \rho_g=0\,1^M
\]
before the orbit lands in \(I_a\) at the point \(\xi_M(a)\). 

We will need the following effective form of this construction. 

\begin{lemma}[Effective landing]
There is a computable map
\[
        \omega\in\{0,1\}^{\mathbb N}\longmapsto a(\omega)\in J
\]
such that
\[
        K(a(\omega))=\rho_f\,\Theta(\omega).
\]
In particular, if \(\omega\) is computable, then \(a(\omega)\) is computable
uniformly from an algorithm computing \(\omega\).
\end{lemma}

\begin{proof}
For each \(a\in J\), let
\[
        \psi_{a,0}:I_a\to L_a,
        \qquad
        \psi_{a,1}:I_a\to R_a
\]
be the two inverse branches of \(g_a\).  For a finite word
\(u=u_0\cdots u_{m-1}\), define
\[
        I_a(u)=
        \psi_{a,u_0}\circ\psi_{a,u_1}\circ\cdots
        \circ\psi_{a,u_{m-1}}(I_a).
\]
Then \(I_a(u)\) is exactly the interval of points whose \(g_a\)-itinerary is \(u\).
The endpoints of \(I_a(u)\) are real-analytic functions of \(a\), obtained by
iterating the analytic inverse branches.  They are computable uniformly in
\(a\) and \(u\), because \(\beta_a,\beta_a',l_a,r_a\) are the unique simple
solutions, in fixed rational isolating intervals, of explicit polynomial
equations involving iterates of \(f_a\).

For a finite word \(u\), define the parameter cylinder
\[
        J(u)=\{a\in J:\xi_M(a)\in I_a(u)\}.
\]
We note that \(J(u)\) is a compact interval, computable uniformly from
\(u\). 
Indeed, if \(b_{u,-}(a)\) and \(b_{u,+}(a)\) denote the two endpoints of
\(I_a(u)\), then the endpoints of \(J(u)\) are the unique zeros inside $J$ of equations of the form
\[
        \xi_M(a)=b_{u,\pm}(a).
\]
Clearly, the family is nested:
\[
        J(uv)\subseteq J(u).
\]
Therefore, for every infinite sequence
\(\omega\), the nested intersection
\[
        \bigcap_{m\geq1} J(\omega_0\cdots \omega_{m-1})
\]
consists of a single parameter.  Call this parameter \(a(\omega)\).

By construction,
\[
        \xi_M(a(\omega))\in C_{a(\omega)}
\]
and its \(g_{a(\omega)}\)-itinerary is exactly \(\omega\).  Before the landing
time, all parameters in \(J\) have the fixed prefix \(\rho_g=0\,1^M\).  Hence
the \(g\)-itinerary of the critical orbit is
\[
        \rho_g\,\omega.
\]
Passing to the kneading sequence for  \(f_a\) gives
\[
        K(a(\omega))
        =
        \rho_f\,\Theta(\omega),
\]
as claimed. Finally, if \(\omega\) is computable, then the endpoints of
\[
        J(\omega_0\cdots\omega_{m-1})
\]
can be computed uniformly from the first \(m\) symbols of \(\omega\).  Since
the lengths of these intervals decrease exponentially, computing such an
interval for \(m\) large enough gives \(a(\omega)\) to any prescribed
precision.  Thus \(a(\omega)\) is computable uniformly from \(\omega\).
\end{proof}

It remains to choose one polytime computable observer whose language can realize every
recursively enumerable set. Fix \(R\geq3\). Let ${\rm{bin}}(n)=(b_1\cdots b_s)_2$ denote the binary expansion of $n$ and define
\[
        \widehat{\rm{bin}}(n)=\widehat b_1\cdots\widehat b_s,
        \quad\text{ with }\quad
        \widehat 0=10,\quad
        \widehat 1=11.
\]
Now put
\[
        Q(n)=0^R1\,\widehat {\rm{bin}}(n)\,01^R
\]
and let the observer be defined by 
\[
        w_n=\Theta(Q(n)).
\]

\begin{lemma}[Synchronization]
Let \(R\geq3\).

\begin{enumerate}
\item Let \(B\) be any word obtained by concatenating \(Q\)-blocks and spacer
blocks \(1^R\).  If \(Q(n)\) occurs in \(B\), then this occurrence is one of
the displayed \(Q\)-blocks.  In particular, \(Q(n)\) cannot occur inside a
different \(Q(m)\), nor across a concatenation boundary.

\item For every infinite word \(\omega\) and every \(n\),
\[
        w_n\text{ occurs in }\rho_f\Theta(\omega)
        \quad\Longleftrightarrow\quad
        Q(n)\text{ occurs in }\omega.
\]
\end{enumerate}
\end{lemma}

\begin{proof}
We first prove the self-synchronization of the \(Q\)-blocks.  Every \(Q(n)\)
begins with the marker \(0^R1\).  Since \(R\geq3\), this marker occurs only at
the beginning of a \(Q\)-block since inside the encoded part \(\widehat{\rm{bin}}(n)\)
all zeros are isolated because
\[
        \widehat 0=10,\qquad \widehat 1=11.
\]
At the junction between \(\widehat{\rm{bin}}(n)\)  and the final delimiter \(01^R\), one
may create two consecutive zeros, but never \(R\) consecutive zeros.  The spacer blocks \(1^R\) contain no zeros.  Hence any occurrence of a word starting with \(0^R1\) must start at the beginning of some \(Q(m)\). Starting at the beginning of \(Q(m)\), the rest of the word is parsed
uniquely.  After the marker \(0^R1\), the encoded part is read in pairs, and
each pair begins with \(1\).  The final delimiter is the first place, at a
pair boundary, where a \(0\) appears instead.  Thus an occurrence of \(Q(n)\)
starting at the beginning of \(Q(m)\) forces
\[
        \widehat n=\widehat m,
        \qquad\text{hence}\qquad
        n=m.
\]
This proves the first claim. For the second claim, recall that
\[
\rho_f=10\,010\,(110)^{M-1}111,
\qquad \Theta(s)=s10.
\]The word \(w_n=\Theta(Q(n))\) begins with \((010)^R\) and contains no \(111\). Since \(\rho_f\) contains no \((010)^R\), an occurrence of \(w_n\) cannot lie entirely in \(\rho_f\). Nor can it cross the junction with \(\Theta(\omega)\): because \(w_n\) starts with \(0\), such an occurrence would have to include the terminal \(111\) of \(\rho_f\). Thus every occurrence lies entirely in \(\Theta(\omega)\).
Moreover, every occurrence must start at a \(\Theta\)-block boundary. Starting at the second position of a block would give initial symbol \(1\); starting at the third position would give third symbol \(1\). Both contradict the initial \(010\) of \(w_n\). Reading the aligned blocks therefore gives an occurrence of \(Q(n)\) in \(\omega\). Conversely, applying \(\Theta\) to any occurrence of \(Q(n)\) in \(\omega\) gives an occurrence of \(w_n\) in \(\rho_f\Theta(\omega)\).
\end{proof}

We now build the realizing sequence.  Since the fixed prefix
\[
        \rho_g=0\,1^M
\]
contains no occurrence of the marker \(0^R1\), no \(Q(n)\) can occur inside
the prefix. Let \(W\subseteq\mathbb N\) be recursively enumerable, and fix an enumeration
machine $M_e$ for \(W\).  Build a computable infinite sequence \(\omega_W\) in
stages.  At stage \(s\), run the enumeration machine for \(s\) steps, list the
numbers which appear for the first time at that stage, in increasing order,
and append \(Q(n)1^R\) for each of them.  If no new number appears at that
stage, append only \(1^R\).  This gives a computable sequence \(\omega_W\),
uniformly from the index $e$. 

By the first part of the synchronization lemma, together with the fact that
\(\rho_g\) contains no occurrence of \(0^R1\),
\[
        Q(n)\text{ occurs in }\omega_W
        \quad\Longleftrightarrow\quad
        n\in W.
\]
Applying the effective landing lemma to \(\omega_W\), we obtain a computable
parameter
\[
        a_W=a(\omega_W)
\]
such that
\[
        K(a_W)=\rho_f\,\Theta(\omega_W).
\]
The synchronization property of \(\Theta\) now gives, for every \(n\),
\[
        w_n\text{ occurs in }K(a_W)
        \quad\Longleftrightarrow\quad
        Q(n)\text{ occurs in }\omega_W
        \quad\Longleftrightarrow\quad
        n\in W.
\]
Therefore
\[
        S_{\mathcal W}(a_W)=W.
\]

%% file: logistic_orbit_adds_ten_tikz.tex
\definecolor{gcurve}{HTML}{7B68EE}
\definecolor{leftfill}{HTML}{F2B66D}
\definecolor{rightfill}{HTML}{6FA8DC}
\definecolor{leftorbit}{HTML}{D97904}
\definecolor{rightorbit}{HTML}{1F77B4}
\definecolor{leftmark}{HTML}{8A5300}
\definecolor{rightmark}{HTML}{1D5F99}
\definecolor{purpleann}{HTML}{5D52BA}

\pgfmathsetmacro{\astar}{3.856800652477765}
\pgfmathsetmacro{\crit}{0.5}

\pgfmathdeclarefunction{fa}{1}{%
  \pgfmathparse{\astar*#1*(1-#1)}%
}
\pgfmathdeclarefunction{ga}{1}{%
  \pgfmathparse{\astar*(\astar*(\astar*#1*(1-#1))*(1-\astar*#1*(1-#1)))*(1-\astar*(\astar*#1*(1-#1))*(1-\astar*#1*(1-#1)))}%
}
\pgfmathsetmacro{\betapcoord}{ga(\crit)}
\pgfmathsetmacro{\betacoord}{ga(\betapcoord)}

\scalebox{0.85}{
\begin{tikzpicture}[
  orbit arrow/.style={-{Stealth[length=1.8mm]}, line width=0.55pt},
  g arrow/.style={-{Stealth[length=1.8mm]}, draw=gcurve, line width=0.65pt, opacity=0.70}
]
\begin{axis}[
  width=12.4cm,
  height=9.4cm,
  xmin=0, xmax=1,
  ymin=0, ymax=1.04,
  axis lines=left,
  xlabel={$x$},
  ylabel={$f_a(x)$},
  xlabel style={yshift=-1.4em},
  ylabel style={rotate=-90, at={(axis description cs:-0.095,1.03)}, anchor=south},
  title={Logistic map, third iterate, and three-step orbit portraits},
  title style={font=\small},
  xtick={0,1},
  xticklabels={$0$,$1$},
  ytick={0,1},
  yticklabels={$0$,$1$},
  tick label style={font=\small},
  label style={font=\small},
  legend style={
    at={(0.98,-0.10)},
    anchor=south east,
    font=\scriptsize,
    draw=gray!35,
    fill=white,
    fill opacity=0.92,
    text opacity=1
  },
  clip=false,
]

\path[fill=leftfill, opacity=0.18]
  (axis cs:\betapcoord,0) rectangle (axis cs:\crit,1.04);
\path[fill=rightfill, opacity=0.17]
  (axis cs:\crit,0) rectangle (axis cs:\betacoord,1.04);

\draw[gray!60, dashed, line width=0.5pt]
  (axis cs:0,0) -- (axis cs:1,1);
\draw[gray!65, dotted, line width=0.55pt]
  (axis cs:\crit,0) -- (axis cs:\crit,1.04);
\draw[gray!65, dotted, line width=0.55pt]
  (axis cs:0,\crit) -- (axis cs:1,\crit);
\draw[leftmark, opacity=0.85, line width=0.45pt]
  (axis cs:\betapcoord,0) -- (axis cs:\betapcoord,1.04);
\draw[rightmark, opacity=0.85, line width=0.45pt]
  (axis cs:\betacoord,0) -- (axis cs:\betacoord,1.04);

\addplot[domain=0:1, samples=260, draw=gcurve, line width=0.75pt, opacity=0.42]
  {ga(x)};
\addlegendentry{$g_a(x)=f_a^3(x)$}

\addplot[domain=0:1, samples=260, draw=black, line width=0.9pt]
  {fa(x)};
\addlegendentry{$f_a(x)=a x(1-x)$}

\addplot[domain=0:1, samples=2, draw=gray!65, dashed, line width=0.5pt]
  {x};
\addlegendentry{$y=x$}

\node[anchor=north, font=\scriptsize, text=leftmark]
  at (axis cs:\betapcoord,-0.055) {$\beta_a'$};
\node[anchor=north, font=\scriptsize, text=gray!55!black]
  at (axis cs:\crit,-0.055) {$c$};
\node[anchor=north, font=\scriptsize, text=rightmark]
  at (axis cs:\betacoord,-0.055) {$\beta_a$};
\node[anchor=east, font=\scriptsize, text=leftmark]
  at (axis cs:-0.035,\betapcoord) {$\beta_a'$};
\node[anchor=east, font=\scriptsize, text=gray!55!black]
  at (axis cs:-0.035,\crit) {$c$};
\node[anchor=east, font=\scriptsize, text=rightmark]
  at (axis cs:-0.035,\betacoord) {$\beta_a$};
\draw[leftmark, opacity=0.55, line width=0.35pt]
  (axis cs:\betapcoord,0) -- (axis cs:\betapcoord,-0.02);
\draw[rightmark, opacity=0.55, line width=0.35pt]
  (axis cs:\betacoord,0) -- (axis cs:\betacoord,-0.02);

\node[font=\small, text=leftmark]
  at (axis cs:0.472549532150440,0.055) {$L_a$};
\node[font=\small, text=rightmark]
  at (axis cs:0.527450467849561,0.055) {$R_a$};

\draw[orbit arrow, draw=leftorbit]
  (axis cs:0.480784672505,0) -- (axis cs:0.480784672505,0.962776121201);
\draw[orbit arrow, draw=leftorbit]
  (axis cs:0.480784672505,0.962776121201) -- (axis cs:0.962776121201,0.962776121201);
\draw[orbit arrow, draw=leftorbit]
  (axis cs:0.962776121201,0.962776121201) -- (axis cs:0.962776121201,0.138221030899);
\draw[orbit arrow, draw=leftorbit]
  (axis cs:0.962776121201,0.138221030899) -- (axis cs:0.138221030899,0.138221030899);
\draw[orbit arrow, draw=leftorbit]
  (axis cs:0.138221030899,0.138221030899) -- (axis cs:0.138221030899,0.459406579806);
\draw[orbit arrow, draw=leftorbit]
  (axis cs:0.138221030899,0.459406579806) -- (axis cs:0.459406579806,0.459406579806);

\foreach \x/\y in {
  0.480784672505/0,
  0.962776121201/0.962776121201,
  0.138221030899/0.138221030899,
  0.459406579806/0.459406579806}
  \addplot[only marks, mark=*, mark size=1.25pt, draw=leftorbit, fill=leftorbit] coordinates {(\x,\y)};

\draw[orbit arrow, draw=rightorbit]
  (axis cs:0.545018767273,0) -- (axis cs:0.545018767273,0.956383626093);
\draw[orbit arrow, draw=rightorbit]
  (axis cs:0.545018767273,0.956383626093) -- (axis cs:0.956383626093,0.956383626093);
\draw[orbit arrow, draw=rightorbit]
  (axis cs:0.956383626093,0.956383626093) -- (axis cs:0.956383626093,0.160882527783);
\draw[orbit arrow, draw=rightorbit]
  (axis cs:0.956383626093,0.160882527783) -- (axis cs:0.160882527783,0.160882527783);
\draw[orbit arrow, draw=rightorbit]
  (axis cs:0.160882527783,0.160882527783) -- (axis cs:0.160882527783,0.520665542740);
\draw[orbit arrow, draw=rightorbit]
  (axis cs:0.160882527783,0.520665542740) -- (axis cs:0.520665542740,0.520665542740);

\foreach \x/\y in {
  0.545018767273/0,
  0.956383626093/0.956383626093,
  0.160882527783/0.160882527783,
  0.520665542740/0.520665542740}
  \addplot[only marks, mark=*, mark size=1.25pt, draw=rightorbit, fill=rightorbit] coordinates {(\x,\y)};

\draw[g arrow] (axis cs:\crit,0.02) -- (axis cs:\crit,\betapcoord);
\draw[g arrow] (axis cs:\crit,\betapcoord) -- (axis cs:\betapcoord,\betapcoord);
\draw[g arrow] (axis cs:\betapcoord,\betapcoord) -- (axis cs:\betapcoord,\betacoord);
\draw[g arrow] (axis cs:\betapcoord,\betacoord) -- (axis cs:\betacoord,\betacoord);
\foreach \x/\y in {
  \crit/\betapcoord,
  \betapcoord/\betapcoord,
  \betapcoord/\betacoord,
  \betacoord/\betacoord}
  \addplot[only marks, mark=*, mark size=1.45pt, draw=gcurve, fill=gcurve, opacity=0.82] coordinates {(\x,\y)};

\node[
  anchor=west,
  align=left,
  font=\scriptsize,
  text=purpleann,
  fill=white,
  draw=gray!35,
  rounded corners=1.5pt,
  inner sep=2.2pt,
  fill opacity=0.90,
  text opacity=1
] (critlabel) at (axis cs:0.61,0.40)
  {critical orbit under $g_a$:\\[-1pt]$c\mapsto\beta_a'\mapsto\beta_a$};
\draw[-{Stealth[length=1.7mm]}, draw=gcurve, line width=0.45pt, opacity=0.75]
  (critlabel.west) -- (axis cs:\betapcoord,\betapcoord);

\node[anchor=west, font=\scriptsize, text=purpleann]
  (gbetap) at (axis cs:0.08,0.69) {$g_a(\beta_a')=\beta_a$};
\draw[-{Stealth[length=1.7mm]}, draw=gcurve, line width=0.5pt, opacity=0.75]
  (gbetap.east) -- (axis cs:\betapcoord,\betacoord);

\node[anchor=west, font=\scriptsize, text=purpleann]
  (gbeta) at (axis cs:0.65,0.895) {$g_a(\beta_a)=\beta_a$};
\draw[-{Stealth[length=1.7mm]}, draw=gcurve, line width=0.5pt, opacity=0.75]
  (gbeta.south west) -- (axis cs:\betacoord,\betacoord);

\node[
  anchor=north west,
  font=\scriptsize,
  fill=white,
  draw=gray!35,
  rounded corners=1.5pt,
  inner xsep=5pt,
  inner ysep=5pt
] at (rel axis cs:0.03,0.965) {$a_\ast\approx 3.856800652478$};

\end{axis}
\end{tikzpicture}
}

%% file: logistic_prefix_landing_tikz.tex
\definecolor{orbitpurple}{HTML}{6F5BD6}
\definecolor{leftfill}{HTML}{F2B66D}
\definecolor{rightfill}{HTML}{6FA8DC}
\definecolor{leftmark}{HTML}{D97904}
\definecolor{rightmark}{HTML}{1F77B4}
\definecolor{darkpurple}{HTML}{4F46AD}

\pgfmathsetmacro{\betapcoord}{0.445097046922302}
\pgfmathsetmacro{\crit}{0.5}
\pgfmathsetmacro{\betacoord}{0.554902953077698}
\pgfmathsetmacro{\midcoord}{0.5*(\betapcoord+\betacoord)}
\pgfmathsetmacro{\symbolheight}{0.812}

\scalebox{0.85}{
\begin{tikzpicture}
\begin{axis}[
  width=12.6cm,
  height=7.1cm,
  xmin=-0.35, xmax=10.8,
  ymin=0.425, ymax=0.855,
  axis lines=left,
  xlabel={$g_a$-time $m$},
  ylabel={$g_a^m(c)$},
  title={Critical orbit producing $\rho_g=0\,1^8$},
  title style={font=\small},
  tick label style={font=\small},
  label style={font=\small},
  xtick={0,1,2,3,4,5,6,7,8,9,10},
  clip=false,
]

\path[fill=gray!60, opacity=0.24]
  (axis cs:-0.35,\betapcoord) rectangle (axis cs:10.8,\betacoord);
\path[fill=leftfill, opacity=0.15]
  (axis cs:-0.35,\betapcoord) rectangle (axis cs:10.8,\crit);
\path[fill=rightfill, opacity=0.14]
  (axis cs:-0.35,\crit) rectangle (axis cs:10.8,\betacoord);

\draw[leftmark, opacity=0.75, line width=0.45pt]
  (axis cs:-0.35,\betapcoord) -- (axis cs:10.8,\betapcoord);
\draw[gray!55, dotted, line width=0.55pt]
  (axis cs:-0.35,\crit) -- (axis cs:10.8,\crit);
\draw[rightmark, opacity=0.75, line width=0.45pt]
  (axis cs:-0.35,\betacoord) -- (axis cs:10.8,\betacoord);

\node[anchor=east, font=\scriptsize, text=leftmark]
  at (axis cs:-0.25,\betapcoord) {$\beta_a'$};
\node[anchor=east, font=\scriptsize, text=gray!55!black]
  at (axis cs:-0.25,\crit) {$c$};
\node[anchor=east, font=\scriptsize, text=rightmark]
  at (axis cs:-0.25,\betacoord) {$\beta_a$};
\node[anchor=west, font=\small]
  at (axis cs:10.62,\midcoord) {$I_a$};

\addplot[
  draw=orbitpurple,
  line width=0.85pt,
  mark=*,
  mark size=1.55pt,
  mark options={draw=orbitpurple, fill=orbitpurple}
] coordinates {
  (0,0.500000000000)
  (1,0.445092408343)
  (2,0.554920186152)
  (3,0.554966981887)
  (4,0.555140920456)
  (5,0.555788352152)
  (6,0.558210685394)
  (7,0.567443778782)
  (8,0.604849887513)
  (9,0.774144656300)
  (10,0.499878620854)
};

\node[
  anchor=center,
  font=\scriptsize,
  text=leftmark,
  fill=white,
  draw=leftmark,
  rounded corners=1.5pt,
  inner xsep=2.8pt,
  inner ysep=1.8pt
] at (axis cs:1,\symbolheight) {$0$};

\node[
  anchor=center,
  font=\scriptsize,
  text=rightmark,
  fill=white,
  draw=rightmark,
  rounded corners=1.5pt,
  inner xsep=2.8pt,
  inner ysep=1.8pt
] at (axis cs:2,\symbolheight) {$1$};
\node[
  anchor=center,
  font=\scriptsize,
  text=rightmark,
  fill=white,
  draw=rightmark,
  rounded corners=1.5pt,
  inner xsep=2.8pt,
  inner ysep=1.8pt
] at (axis cs:3,\symbolheight) {$1$};
\node[
  anchor=center,
  font=\scriptsize,
  text=rightmark,
  fill=white,
  draw=rightmark,
  rounded corners=1.5pt,
  inner xsep=2.8pt,
  inner ysep=1.8pt
] at (axis cs:4,\symbolheight) {$1$};
\node[
  anchor=center,
  font=\scriptsize,
  text=rightmark,
  fill=white,
  draw=rightmark,
  rounded corners=1.5pt,
  inner xsep=2.8pt,
  inner ysep=1.8pt
] at (axis cs:5,\symbolheight) {$1$};
\node[
  anchor=center,
  font=\scriptsize,
  text=rightmark,
  fill=white,
  draw=rightmark,
  rounded corners=1.5pt,
  inner xsep=2.8pt,
  inner ysep=1.8pt
] at (axis cs:6,\symbolheight) {$1$};
\node[
  anchor=center,
  font=\scriptsize,
  text=rightmark,
  fill=white,
  draw=rightmark,
  rounded corners=1.5pt,
  inner xsep=2.8pt,
  inner ysep=1.8pt
] at (axis cs:7,\symbolheight) {$1$};
\node[
  anchor=center,
  font=\scriptsize,
  text=rightmark,
  fill=white,
  draw=rightmark,
  rounded corners=1.5pt,
  inner xsep=2.8pt,
  inner ysep=1.8pt
] at (axis cs:8,\symbolheight) {$1$};
\node[
  anchor=center,
  font=\scriptsize,
  text=rightmark,
  fill=white,
  draw=rightmark,
  rounded corners=1.5pt,
  inner xsep=2.8pt,
  inner ysep=1.8pt
] at (axis cs:9,\symbolheight) {$1$};

\node[
  anchor=west,
  align=left,
  font=\scriptsize,
  fill=white,
  draw=gray!35,
  rounded corners=1.5pt,
  inner sep=2.5pt,
  fill opacity=0.92,
  text opacity=1
] (landinglabel) at (axis cs:8.40,0.665)
  {$\xi_M(a)=g_a^{M+2}(c)\in I_a$};
\draw[-{Stealth[length=1.8mm]}, draw=orbitpurple, line width=0.55pt]
  (landinglabel.east) -- (axis cs:10,0.499878620854);

\node[
  anchor=west,
  font=\scriptsize,
  text=darkpurple
] (shadowlabel) at (axis cs:4.35,0.595)
  {orbit shadows the repelling fixed point};
\draw[-{Stealth[length=1.8mm]}, draw=orbitpurple, line width=0.55pt, opacity=0.80]
  (shadowlabel.west) -- (axis cs:3.2,0.5556);

\end{axis}
\end{tikzpicture}
}

%% file: universal_implies_chaotic.tex
\section{Universality implies chaos}

We prove that universal parameters must have positive
topological entropy. We shall use the following standard consequence of kneading theory and the
period-doubling classification for the quadratic family; see
\cite{kneading,deMeloVanStrien1993}.

\begin{proposition}
If \(h_{\mathrm{top}}(f_a)=0\), then one of the following holds:
\begin{enumerate}
\item the critical orbit is attracted to a periodic orbit of period \(2^q\),
for some \(q\geq0\);
\item \(a=a_F\), the Feigenbaum parameter.
\end{enumerate}
\end{proposition}

\begin{theorem}
If \(a\) is universal, then
\[
        h_{\mathrm{top}}(f_a)>0.
\]
\end{theorem}

\begin{proof}
Suppose first that the critical orbit is attracted to a periodic orbit of
period \(2^q\).  In this case $a$ is regular and therefore has a decidable kneading sequence. 

It remains to consider the Feigenbaum parameter.  In this case the kneading
sequence is the period-doubling sequence.  Equivalently, it is the fixed point
of the substitution
\[
        \tau(1)=10,\qquad \tau(0)=11,
\]
namely
\[
        K(a_F)=\lim_{n\to\infty}\tau^n(1)
              =1011101010111011\cdots .
\]
This symbolic description is classical. To see that its language is decidable, let \(u\) be a word of length
\(m\), and choose \(N\) such that \(2^N\geq m\).  Since
\[
        K(a_F)=\tau^N(K(a_F)),
\]
the sequence is a concatenation of blocks \(\tau^N(0)\) and \(\tau^N(1)\),
each of length \(2^N\).  Therefore any occurrence of \(u\) is contained in
\(\tau^N(ab)\) for some adjacent pair \(ab\) occurring in \(K(a_F)\).  The
only adjacent pairs are
\[
        01,\qquad 10,\qquad 11.
\]
Thus \(u\) occurs in \(K(a_F)\) if and only if it occurs in one of the three
finite words
\[
        \tau^N(01),\qquad \tau^N(10),\qquad \tau^N(11).
\]
This gives a finite decision procedure. Since $K(a)$ is decidable in both cases, so is 
 \(S_{\mathcal W}(a)\) for any computable observer \(\mathcal W=(w_n)\). Therefore, no
zero-entropy parameter can be universal.  
\end{proof}